\documentclass[12pt,british]{amsart}

\usepackage[latin1]{inputenc}
\usepackage[T1]{fontenc}
\usepackage{babel,eucal,url,amssymb,stmaryrd,booktabs,%
enumerate,amscd,paralist}

\theoremstyle{plain}
\newtheorem{theorem}{Theorem}[section]
\newtheorem{proposition}[theorem]{Proposition}
\newtheorem{lemma}[theorem]{Lemma}

\theoremstyle{definition}
\newtheorem{definition}[theorem]{Definition}

\theoremstyle{remark}

\newtheorem{example}[theorem]{Example}

\numberwithin{equation}{section}

\DeclareMathOperator{\alg}{alg}

\DeclareMathOperator{\inv}{inv}
\DeclareMathOperator{\C}{C}
\DeclareMathOperator{\id}{id}
\DeclareMathOperator{\NC}{NC}

\newcommand{\norm}[1]{\left\lVert #1\right\rVert}

\begin{document}

\title[Amalgamated R- and S-transform \today]{A Banach algebra approach
to amalgamated R- and S-transforms}

\author{Lars Aagaard}
\address[Aagaard]{Department of Mathematics and Computer Science\\
University of Southern Denmark\\
Campusvej 55\\
DK-5230 Odense M\\
Denmark}
\email{laa@imada.sdu.dk}

%\author{Advisor: Uffe Haagerup}
%\address[Haagerup]{Department of Mathematics and Computer Science\\
%University of Southern Denmark\\
%Campusvej 55\\
%DK-5230 Odense M\\
%Denmark}
%\email{Haagerup@imada.sdu.dk}

\begin{abstract}
We give a Banach algebra approach to the additiveness property of Voiculescu's amalgamated R-transform. We also define an amalgamated S-transform and prove that it is multiplicative on products of amalgamated free random variables when the algebra of amalgamation is commutative.
\end{abstract}

%\subjclass[2000]{Primary 53C10; Secondary 17B10, 53C25, 53C29}

\maketitle

% \newpage
% \tableofcontents
% newpage

\section{Introduction}

If $A$ is a unital algebra and $\phi$ is a unital linear functional on $A$ we will say that $(A,\phi)$ is a free probability space. The elements of $A$ are then the free random variables, and $\phi$ serves as the ``probability measure'' on $A$. Given two random variables $a_1, a_2\in A$ we say that $a_1$ and $a_2$ are free wrt. $\phi$ if for all $n\in \mathbb N$ and all polynomials $(p_i)_{i=1}^n$ we have
\begin{equation*}
  \phi(p_1(a_{i_1})p_2(a_{i_2})\cdots p_n(a_{i_n}))=0,
\end{equation*}
whenever $i_1\neq i_2 \neq \cdots i_n\neq i_1$ and $\phi(p_j(a_j))=0$ for $j=1,\ldots n$. 

It turns out that freeness of $a_1$ and $a_2$ and knowledge of $\phi(a_1^k)$ and $\phi(a_2^k)$ for all $k\in \mathbb N$ is enough to compute the distributions of $a_1+a_2$ and $a_1a_2$, i.e. finding the moments $\phi((a_1+a_2)^k)$ and $\phi((a_1a_2)^k)$ respectively for all $k\in \mathbb N$. This is done in turn by use of Voiculescu's $R$- and $S$-transforms. 

For fixed $a\in A$ one define the Cauchy-tranform of $a$ by
\begin{equation*}
  G_a(z)= \phi((z-a)^{-1})
\end{equation*}
for $z\in \mathbb C$ whenever this makes sense, and also define
\begin{equation*}
  \psi_a(z) = \phi((1-za)^{-1})-1
\end{equation*}
for $z\in \mathbb C$. The $R$- and $S$-transforms are then defined by
\begin{equation*}
  R_a(z) = G^{\langle -1 \rangle}_a(z) - z^{-1}
\end{equation*} 
and 
\begin{equation*}
  S_a(z) = \frac{1+z}{z}\psi^{\langle -1 \rangle}_a(z)
\end{equation*}
for $z \in \mathbb C$, where superscript $\langle -1\rangle$ denotes inversion wrt. composition.
Voiculescu then proved the following formulas \cite{Voi1}, \cite{Voi2}:
\begin{equation} \label{Rtransformen}
  R_{a_1+a_2}(z)=R_{a_1}(z)+R_{a_2}(z) 
\end{equation}
and 
\begin{equation} \label{Stransformen}
  S_{a_1\cdot a_2}(z)=S_{a_1}(z)\cdot S_{a_2}(z).
\end{equation}
when $a_1$ and $a_2$ are free wrt. $\phi$ and $z\in \mathbb C$ is chosen suitably.

Alternative proofs of the above formulas can be found in \cite{Haa}, \cite{Sp1} and \cite{NS1}. A nice combinatorial review can be found in \cite{Combinatorics}.

Instead of considering a ``classical'' free probability space one can consider the more general notion of an amalgamated free probability space, i.e. let $\mathcal A$ be a unital Banach algebra, 
$1\in \mathcal
B\subset \mathcal A$ a unital Banach sub-algebra of $\mathcal A$, and 
$E:\mathcal A\to \mathcal B$ a conditional expectation. Thus 
$E$ is linear, $E(b)=b$
for all $b\in \mathcal B$, $E$ is norm-decreasing, and $E$ has the
$\mathcal B$-bimodule property;  
\begin{equation*}
  E(b_1ab_2)=b_1E(a)b_2
\end{equation*}
for $b_1, b_2\in \mathcal B$ and $a\in \mathcal A$. We will say that $(\mathcal B \subset \mathcal A,E)$ is a $\mathcal B$-probability space.

By replacing freeness above by the corresponding freeness with amalgamation Voiculescu proved an amalgamated version of (\ref{Rtransformen}) in \cite{Voi3}, and Speicher reproved this by combinatorial means in \cite{AMS}. The definition of freeness wih amalgamation is as follows.
Let $a_1, a_2 \in \mathcal A$ be random variables in $\mathcal A$. We will say that 
$a_1$ and $a_2$ are free with amalgamation over $\mathcal B$ wrt. $E$ (or free wrt. $\mathcal B$ or simply  
$\mathcal B$-free) if for all $n\in\mathbb N$ and for all 
$x_j\in \alg (\mathcal B,a_{i_j})$ we have
\begin{equation} \label{eq:1.3}
  E(x_1x_2 \cdots x_n)=0,
\end{equation}
whenever $i_1\neq i_2 \neq \cdots \neq i_n$ and $E(x_j)=0$ for all $j=1,\ldots,n$. In other words the indices has to be alternating since we are only dealing with two random variables.

Note, that if the condition $E(x_j)=0$, $1\leq j\leq n$ is relaxed to $E(x_j)=0$ for $2\leq j \leq n-1$, (\ref{eq:1.3}) should be changed to 
\begin{equation} \label{eq:1.4}
  E(x_1 x_2 \cdots x_{n-1})=
  \begin{cases}
    E(x_1)E(x_{n-1}), & n=2 \\
    0, & n\geq 3
  \end{cases}.
\end{equation}
This can be seen by writing $x_1=x_1^0+ E(x_1)$ and $x_n=x_n^0+E(x_n)$, where $E(x^0_1)=E(x^0_n)=0$.

In the unital Banach-algebra seting we give a new proof of the additiveness property of the amalgamated $R$-transform. We also prove an amalgamated version of (\ref{Stransformen}) for the amalgamated $S$-transform in the case where $\mathcal B$ is abelian.
Our methods are strongly inspired by section 3 of \cite{Haa} and as in \cite{Haa} we give concrete neighboorhoods where amalgamated versions of (\ref{Rtransformen}) and (\ref{Stransformen}) are valid.

\section{Amalgamated R-transform in Banach-algebras}

Throughout this paper we let $\mathcal A$ be a unital Banach algebra,  
$1\in \mathcal
B\subset \mathcal A$ a unital Banach sub-algebra of $\mathcal A$, and 
$E:\mathcal A\to \mathcal B$ a conditional expectation, so that $(\mathcal B\subset \mathcal A,E)$ is a $\mathcal B$-probability space.  

For $\epsilon>0$ we denote by $\mathcal
B(0,\epsilon)$ the $\epsilon$-neighborhood of $0$ in $\mathcal B$ and by
$\mathcal 
B(0,\epsilon)_{\inv}$ the invertible
elements in $\mathcal B$ of norm strictly less than
$\epsilon$. Finally $\overline{\mathcal B}(0,\epsilon)$ will denote the
norm-closure of $\mathcal B(0,\epsilon)$.

Let $a\in \mathcal A$ be fixed. Define
the function $g_a:\mathcal B(0,\tfrac{1}{\norm{a}})\to \mathcal B$ by 
\begin{equation}
  \label{gdef}
  g_a(b) :=  bE((1-ab)^{-1}) =E((1-ba)^{-1})b
\end{equation}
for $b\in \mathcal B(0,\tfrac{1}{\norm{a}})$. For 
$b\in \mathcal B(0,\tfrac{1}{\norm{a}})$ observe that by the Carl
Neumann series  $g_a(b)$ is just the
absolute convergent sum 
\begin{equation*} 
  g_a(b) = b + bE(a)b + bE(aba)b + \cdots.
\end{equation*}
Of course, if $b\in \mathcal B(0,\tfrac{1}{\norm{a}})_{\inv}$ then $g_a(b) = E((b^{-1}-a)^{-1})$ is
nothing but $G_a(b^{-1})$, where $G_a$ is the amalgamated 
Cauchy-transform of $a$.

\begin{lemma} \label{inject}
Let $g_a:\mathcal B(0,\tfrac{1}{\norm{a}}) \to \mathcal B$ be the
function defined by
(\ref{gdef}). Then $g_a$ is 1-1 on $\mathcal
B(0,\tfrac{1}{4\norm{a}})$. 
\end{lemma}

\begin{proof}
Let $b_1,b_2\in \mathcal B$ such that $\norm{b_i} <
\tfrac{\alpha}{\norm{a}}$ for $0<\alpha<1$ to be determined. The idea
is to determine $\alpha$ such that
\begin{equation} \label{gadiff}
\norm{\bigl((g_a(b_1)-g_a(b_2)\bigr) -\bigl(b_1-b_2\bigr)} < \norm{b_1-b_2}
\end{equation}
for $\norm{b_1},\norm{b_2}< \tfrac{\alpha}{\norm{a}}$.

Observe that
\begin{multline*}
 g_a(b_1) -g_a(b_2)  = (b_1-b_2) + E(b_1ab_1)-E(b_2ab_2) \\ 
   +
 E(b_1ab_1ab_1)- E(b_2ab_2ab_2) + \cdots \\
  =  (b_1-b_2) + E((b_1-b_2)ab_1) + E(b_2a(b_1-b_2)) \\
     + E(b_1-b_2)ab_1ab_1) + E(b_2a(b_1-b_2)ab_2)  \\ +
   E(b_2ab_2a(b_1-b_2)) + \cdots,
\end{multline*}
so we can estimate the left-hand side of (\ref{gadiff}) by
\begin{multline*}
\norm{g_a(b_1)-g_a(b_2) -(b_1-b_2)}  \leq 
2\norm{a}\max\{\norm{b_1},\norm{b_2}\}\norm{b_1-b_2} \\
  + 3 \norm{a}^2\max\{\norm{b_1},\norm{b_2}\}^2\norm{b_1-b_2} 
  + \cdots   \\
 \leq  \norm{b_1-b_2} \left(\sum_{n=0}^\infty(n+1)\alpha^n -1\right)
\\
 =  \left(\frac{1}{(1-\alpha)^2}-1\right)\norm{b_1-b_2}. 
\end{multline*}
We infer that if $\alpha< \tfrac{1}{4}$ then
(\ref{gadiff}) is satisfied.
\end{proof}

We are interested in a neighborhood $\mathcal
B(0,\tfrac{\beta}{\norm{a}})$, where $0<\beta<1$ is to be determined, such
that $g_a$ maps $\mathcal
B(0,\tfrac{1}{4\norm{a}})$ onto a neighborhood of $0$ which contains
the neighborhood $\mathcal B(0,\tfrac{\beta}{\norm{a}})$. For this we
need the following lemma that enables us to use the Inverse Function
Theorem \cite[Theorem 3.6.3]{VLH} on $g_a$.

\begin{lemma} \label{differentiale} Let $a\in \mathcal A$ be fixed.
  Define $g_a:\mathcal B(0,\tfrac{1}{\norm{a}})\to \mathcal B$ by
  $g_a(b) = bE((1-ab)^{-1})$ for $b\in \mathcal 
  B(0,\tfrac{1}{\norm{a}})$. Then $g_a$ is Fr\'echet differentiable, and
  the differential of $g_a$ is 
\begin{equation*}
  Dg_a(b):h\mapsto E((1-ba)^{-1}h(1-ab)^{-1})
\end{equation*}
for $b\in \mathcal B(0,\tfrac{1}{\norm{a}})$ and $h\in \mathcal B$.
Furthermore the differential, $Dg_a:\mathcal B(0,\tfrac{1}{\norm{a}})\to 
\mathcal B$, is continuous, so that $g_a$ is differentiable of class 
$\C^1$. If $b \in \mathcal B(0,\tfrac{\alpha}{\norm{a}})$ for 
$0<\alpha<1$ then
\begin{equation}
  \label{Dgavurdering}
  \norm{Dg_a(b)-Dg_a(0)} < \frac{\alpha(2-\alpha)}{(1-\alpha)^2}. 
\end{equation}
\end{lemma}
 
\begin{proof}
If $b_1, b_2\in \mathcal B(0,\tfrac{1}{\norm{a}})$ then 
\begin{multline}
  \norm{(1-b_1a)^{-1} - (1-b_2a)^{-1}}   \\
= \norm{(1-b_1a)^{-1}\bigl((1-b_2a)-(1-b_1a)\bigr)(1-b_2a)^{-1}}  \\
= \norm{(1-b_1a)^{-1}(b_2-b_1)a(1-b_2a)^{-1}} \label{diff1} \\
\leq \norm{(1-b_1a)^{-1}} \norm{(1-b_2a)^{-1}}\norm{a} \norm{b_2-b_1} \\
\leq \frac{1}{1-\norm{a}\norm{b_1}} \frac{1}{1-\norm{b_2}\norm{a}}\norm{a} \norm{b_2-b_1} \to 0 
\end{multline}
as $b_2\to b_1$ in norm.
Also 
\begin{multline} \label{diff2}
  g_a(b_1)-g_a(b_2)  =  E((1-b_1a)^{-1}b_1 - b_2(1-ab_2)^{-1})  \\
    =  E\left((1-b_1a)^{-1}\bigl(b_1(1-ab_2)-(1-b_1a)b_2\bigr)(1-ab_2)^{-1}\right)
  \\   =  E((1-b_1a)^{-1}(b_1-b_2)(1-ab_2)^{-1}). 
\end{multline}
Combining (\ref{diff1}) and (\ref{diff2}) we conclude that for $b\in
\mathcal B(0,\tfrac{1}{\norm{a}})$ fixed, and $h\in \mathcal B$  of small norm the first part of the lemma now follows since 
\begin{multline*}
  g_a(b+h)-g_a(b)  =  E((1-(b+h)a)^{-1}h(1-ab)^{-1}) \\
 =  E((1-ba)^{-1}h(1-ab)^{-1}) \\ 
 + E\left(\left( (1-(b+h)a)^{-1}- (1-ba)^{-1}\right) h (1-ab)^{-1}\right) \\
 =   E((1-ba)^{-1}h(1-ab)^{-1}) + O(\norm{h}^2).
\end{multline*}

To see that $Dg_a:b\mapsto Dg_a(b)$ is continuous for $b\in \mathcal B(0,\tfrac{1}{\norm{a}})$ we observe that if $b_1,b_2\in \mathcal B(0,\tfrac{1}{\norm{a}})$ then (\ref{diff1}) implies
\begin{multline} \label{Dgavur}
  \norm{Dg_a(b_1)-Dg_a(b_2)} 
 =  \sup_{\norm{h}\leq 1} \norm{Dg_a(b_1)(h)-Dg_a(b_2)(h)} \\
 =  \sup_{\norm{h}\leq 1} \norm{E((1-b_1a)^{-1}h(1-ab_1)^{-1})-
E((1-b_2a)^{-1}h(1-ab_2)^{-1}))} \\
 \leq \sup_{\norm{h}\leq 1} \norm{(1-b_1a)^{-1}h(1-ab_1)^{-1}-
(1-b_1a)^{-1}h(1-ab_2)^{-1})} \\
 + \sup_{\norm{h}\leq 1} \norm{(1-b_1a)^{-1}h(1-ab_2)^{-1}-
(1-b_2a)^{-1}h(1-ab_2)^{-1})} \\
\leq \norm{(1-b_1a)^{-1}} \norm{(1-ab_1)^{-1} - (1-ab_2)^{-1}} \\
 +   \norm{(1-b_1a)^{-1} - (1-b_2a)^{-1}} \norm{(1-ab_2)^{-1}} \to 0 
\end{multline}
as $b_2\to b_1$ in norm.

Let $0<\alpha<1$. For $b\in \mathcal B(0,\tfrac{\alpha}{\norm{a}})$ 
we use (\ref{Dgavur}) (with $b_1=b$ and $b_2=0$) to see that
\begin{multline*}
 \norm{Dg_a(b)-Dg_a(0)}  \leq  
 \norm{(1-ba)^{-1}} \norm{(1-ab)^{-1} - 1} \\
 +   \norm{(1-ba)^{-1} - 1}  \\
\leq \left(\frac{1}{1-\norm{b}\norm{a}}+1\right)
 \left(\sum_{n=1}^\infty \norm{a}^n\norm{b}^n\right) \\
  =
 \left(\frac{1}{1-\alpha}+1\right)\frac{\alpha}{1-\alpha}  
  =  \frac{\alpha(2-\alpha)}{(1-\alpha)^2}.
\end{multline*}
\end{proof}

\begin{proposition} \label{domain}
Let $\mathcal A$ be a unital Banach algebra, $1\in \mathcal
B\subset \mathcal A$ a unital Banach sub-algebra, and $E:\mathcal A\to
\mathcal B$ a
conditional expectation. Let $a\in \mathcal A$ be a fixed element. The
function $g_a:\mathcal 
B(0,\tfrac{1}{\norm{a}})\to \mathcal B$ defined by
$g_a(b)=bE((1-ab)^{-1})$ for $b\in \mathcal B(0,
\tfrac{1}{\norm{a}})$ is 
a bijection of the neighborhood $\mathcal B(0,\tfrac{1}{4\norm{a}})$
onto a neighborhood of $0$ which contains $\mathcal
B(0,\tfrac{1}{11\norm{a}})$. Furthermore  
\begin{gather}
%B(0,\tfrac{1}{14\norm{a}})_{\inv}  \subseteq  g_a(B(0,\tfrac{1}{7\norm{a}})_{\inv} )   \\
% g_a\left(\mathcal B(0,\tfrac{1}{9\norm{a}})_{\inv} \right)   
%\subseteq \mathcal B(0,\tfrac{1}{8\norm{a}})_{\inv} \\
 g_a^{<-1>}\left(\mathcal B(0,\tfrac{1}{11\norm{a}}) \right)
   \subseteq \mathcal B(0,\tfrac{2}{11\norm{a}}) \label{2.7} \\
 g_a^{<-1>}\left(\mathcal B(0,\tfrac{1}{11\norm{a}})_{\inv} \right)
   \subseteq \mathcal B(0,\tfrac{2}{11\norm{a}})_{\inv}
   \label{neighboor}
\end{gather}
\end{proposition}

The proof is actually just an inspection of some of the
proof of the Inverse Function Theorem, and can be found in any
standard text book on the subject. The proof we inspect here is taken from \cite{VLH}. Since the estimates are of importance to us 
we include a proof.

\begin{proof}
Define $T:\mathcal B(0,\tfrac{1}{\norm{a}})\to \mathcal B$ by $T(b) =
b- g_a(b)$ for $b\in B(0,\tfrac{1}{\norm{a}})$, and observe that
$T(0)=0$. By lemma \ref{differentiale} we have
\begin{equation*}
DT(0) = D\id_\mathcal B(0) - Dg_a(0) = D\id_\mathcal B(0)-D\id_\mathcal
B(0)=0. 
\end{equation*}
By lemma \ref{differentiale} $Dg_a$ is continuous so $DT$ is
also continuous, and
\begin{equation*}
\norm{DT(b)}  =  \norm{Dg_a(b) - Dg_a(0)} <
\frac{\alpha(2-\alpha)}{(1-\alpha)^2} 
\end{equation*} 
for $b\in \mathcal B(0,\tfrac{\alpha}{\norm{a}})$ and
$0<\alpha<1$. If we choose $\alpha= \tfrac {2+\epsilon}{11}$ for
$\epsilon>0$ small enough we have
\begin{equation*}
\norm{DT(b)} <\frac{\tfrac{2+\epsilon}{11}(2-\tfrac{2+\epsilon}{11}
)}{(1-\tfrac{2+\epsilon}{11})^2} < \frac 1 2. 
\end{equation*}
By the Mean Value Theorem \cite{VLH} we conclude that
\begin{equation*}
\norm{T(b)}\leq \frac 1 2 \norm{b}
\end{equation*}
for $b\in \mathcal B(0,\tfrac{2+\epsilon}{11\norm{a}})$.

Now let $\tilde b\in \mathcal \overline{\mathcal
B}(0,\tfrac{2+\epsilon}{22\norm{a}})$ and define 
\begin{equation*}
T_{\tilde{b}}(b)= \tilde{b} - T(b) = \tilde{b}+b-g_a(b)
\end{equation*}
for $b\in \overline{\mathcal B}(0,\tfrac{2+\epsilon}{11\norm{a}})$. For $b\in
\overline{\mathcal B}(0,\tfrac{2+\epsilon}{11\norm{a}})$ we have
\begin{equation*}
\norm{T_{\tilde{b}}(b)} = \lVert \tilde b - T(b)\rVert  \leq \lVert\tilde b\rVert +
\norm{T(b)} \leq \tfrac{2+\epsilon}{11\norm{a}}, 
\end{equation*}
so actually $T_{\tilde{b}}$ is a continuous map from the Banach space
$\overline{\mathcal B}(0,\tfrac{2+\epsilon}{11\norm{a}})$ into itself. Also
$T_{\tilde{b}}:\overline{\mathcal B}(0,\tfrac{2+\epsilon}{11\norm{a}})\to
\overline{\mathcal B}(0,\tfrac{2+\epsilon}{11\norm{a}})$ is a contraction since
for $b_1,b_2\in 
\overline{\mathcal B}(0,\tfrac{2+\epsilon}{11\norm{a}})$ we have
\begin{equation*}
\norm{T_{\tilde{b}}(b_1)-T_{\tilde{b}}(b_2)} = \norm{T(b_1)-T(b_2)} \leq
\tfrac 1 2 \norm{b_1-b_2},
\end{equation*}
again by the Mean Value Theorem \cite{VLH}. Banach's Fix-point
Theorem \cite{VLH} now implies that there exists a unique fix-point $b\in
\overline{\mathcal B}(0,\tfrac{2+\epsilon}{11\norm{a}})$ of
$T_{\tilde{b}}$. Thus $g_a(b)=\tilde{b}$, so we have shown that
\begin{equation*}
\overline{\mathcal B}(0,\tfrac{2+\epsilon}{22\norm{a}}) \subseteq
g_a(\overline{\mathcal B}(0,\tfrac{2+\epsilon}{11\norm{a}})).
\end{equation*}
Of course an argument 
similar to the above also shows that 
\begin{equation} \label{mindre}
\overline{\mathcal
B}(0,\tfrac{1}{11\norm{a}})\subseteq g_a(\overline{\mathcal
B}(0,\tfrac{2}{11\norm{a}})).
\end{equation} 
Now assume that $\norm{b} =
\tfrac{2}{11\norm{a}}$. If $\norm{g_a(b)} < \tfrac{1}{11\norm{a}}$ then
since $g_a$ is continuous we could also find a $b_1\in \mathcal B$ such
that $\tfrac{2}{11\norm{a}}<\norm{b_1}< \tfrac{2+\epsilon}{11\norm{a}}$ and
such that $\norm{g_a(b_1)} < \tfrac{1}{11\norm{a}}$. But by
(\ref{mindre}) we could also find a $b_2\in \overline{\mathcal
B}(0,\tfrac{2}{11\norm{a}})$ such that $g_a(b_2) = g_a(b_1)$ thus
contradicting the injectivity of $g_a$ on $\mathcal
B(0,\tfrac{1}{4\norm{a}})$ that follows from lemma
\ref{inject}. We have now shown that $g_a$ is a bijection of the
neighborhood $\mathcal B(0,\tfrac{2}{11\norm{a}})$ 
onto a neighborhood of $0$ which contains $\mathcal
B(0,\tfrac{1}{11\norm{a}})$ and this is exactly (\ref{2.7}). 

To prove (\ref{neighboor}), define $\Psi_a(b) := \sum_{n=1}^\infty
E((ab)^n)$ for $b\in \mathcal
B(0,\tfrac{1}{\norm{a}})$. For $b\in \mathcal
B(0,\tfrac{1}{2\norm{a}})$ we have
\begin{equation*}
\norm{\Psi_a(b)} \leq \sum_{n=1}^\infty (\norm{a}\norm{b})^n <
\sum_{n=1}^\infty \frac{1}{2^n} =1,
\end{equation*}
so $1+ \Psi_a(b)$ is invertible. Thus if $b\in \mathcal
B(0,\tfrac{1}{2\norm{a}})$ we have
\begin{equation*}
g_a(b) = bE((1-ab)^{-1}) = b(1+\Psi_a(b)),
\end{equation*}
so invertability of $g_a(b)$ is equivalent to invertability of $b$ and
thus (\ref{neighboor}) follows.
\end{proof}

Proposition \ref{domain} now assures well-definedness when we to define Voiculescu's amalgamated R-transform in the Banach algebra setting. 
\begin{definition} \label{Rdef_inv}
  Let $\mathcal A$ be a unital Banach algebra, $1\in \mathcal
  B\subset \mathcal A$ a unital Banach sub-algebra, and $E:\mathcal
  A\to \mathcal B$ a conditional expectation. Let $a\in \mathcal 
  A$ be a fixed non-zero element. Then the amalgamated $R$-transform of $a$ is defined 
by
  \begin{equation}
    \label{R-transform}
    R_a(b) := \left(g_a^{<-1>}(b)\right)^{-1} - b^{-1}
  \end{equation}
for $b\in \mathcal B(0, \tfrac{1}{11\norm{a}})_{\inv}$.
\end{definition}
Again we can also use proposition \ref{domain} to show the additiveness of the 
amalgamated $R$-transform on free random variables. 
\begin{theorem}
  Let $\mathcal A$ be a unital Banach algebra, $1\in \mathcal
  B\subset \mathcal A$ a unital Banach sub-algebra, and $E:\mathcal
  A\to \mathcal B$ a
  conditional expectation. Let $a_1, a_2\in \mathcal A$ be two $\mathcal
  B$-free random variables with respect to $E$. Let $b\in 
  \mathcal
  B(0,\min(\tfrac{1}{11\norm{a_1+a_2}},\tfrac{1}{11\norm{a_1}},\tfrac{1}{11\norm{a_2}}))_{\inv}$.  
  Then the amalgamated $R$-transform satisfies
  \begin{equation*}
    R_{a_1+a_2}(b) = R_{a_1}(b) + R_{a_2}(b).
  \end{equation*}
\end{theorem} 

\begin{proof}
 In the following proof $i$ always denotes an index such that $i\in \{1,2\}$.
  Recall that $g^{<-1>}_{a_i}$ is defined and one-to-one in
  $B(0,\tfrac{1}{11\norm{a_i}})_{\inv}$ by proposition
  \ref{domain}. Let $b\in
  \mathcal B(0,\min(\tfrac{1}{11\norm{a_1+a_2}},\tfrac{1}{11\norm{a_1}},\tfrac{1}{11\norm{a_2}}))_{\inv}$ 
  and  
let $b_1,b_2\in \mathcal B$ be the uniquely determined elements that satisfy
 $b_i\ \in \mathcal B(0,\tfrac{2}{11\norm{a_i}})_{\inv}$ and 
  \begin{equation} \label{gaib=0}
    b= g_{a_i}(b_i),
  \end{equation}
for $i\in \{1,2\}$.
Observe that
\begin{equation*}
b_i^{-1}b = b_i^{-1} b_i E((1-a_ib_i)^{-1}) =
1+\sum_{n=1}^\infty E((a_ib_i)^n),
\end{equation*}
and thus 
\begin{equation*}
b_1^{-1}b+b_2^{-1}b -1 = 1 + \sum_{n=1}^\infty E((a_1b_1)^n) +
\sum_{n=1}^\infty E((a_2b_2)^n), 
\end{equation*} so we infer that $b_1^{-1}b+b_2^{-1}b -1$ is invertible
and that
\begin{equation*}
\norm{(b_1^{-1}b+b_2^{-1}b -1)^{-1}} \leq 1+2\frac{\frac{2}{11}}{1-\frac{2}{11}} <2,
\end{equation*}
because $b_i\in \mathcal B(0,\tfrac{2}{11\norm{a_i}})$.
If we define $b_3= b(b_1^{-1}b + b_2^{-1}b -1)^{-1}$ then obviously
$b_3\in \mathcal B(0,\tfrac{2}{11\norm{a_1+a_2}})_{\inv}$ because $b\in\mathcal
B(0,\tfrac{1}{11\norm{a_1+a_2}})_{\inv}$. Furthermore 
\begin{equation} \label{b3def}
b_3^{-1}= 
b_1^{-1}+b_2^{-1}-b^{-1}.
\end{equation}
Define $A_1(b_1) = (b_1^{-1}-a_1)^{-1} -b$ and 
$A_1(b_2) = (b_2^{-1}-a_2)^{-1} -b$. By (\ref{gaib=0}) we have
$E(A_1(b_1)) =E(A_2(b_2))=0$. Note that
\begin{multline*}
  (b_1^{-1}-a_1)b\left(1-b^{-1}A_1(b_1)b^{-1}A_2(b_2)\right)(b_2^{-1}-a_2) \\ 
 = (b_1^{-1}-a_1)\left(b-((b_1^{-1}-a_1)^{-1}-b)b^{-1}
 ((b_2^{-1}-a_2)^{-1}-b)\right)(b_2^{-1}-a_2) 
 \\ 
=  (b_1^{-1}-a_1)\Bigl(b- (b_1^{-1}-a_1)^{-1}b^{-1}(b_2^{-1}-a_2)^{-1}
      \\ +(b_1^{-1}-a_1)^{-1} +(b_2^{-1}-a_2)^{-1} -b\Bigr)(b_2^{-1}-a_2) \\
  = -b^{-1} +(b_1^{-1}-a_1) +(b_2^{-1}-a_2) \\
  =  b_3^{-1} -( a_1+a_2). ~~~~~~~~~~~~~~~~~~~~~~~~~~~
\end{multline*}
Inverting we have
\begin{multline}
\label{inverti}
 (b_3^{-1} -( a_1+a_2))^{-1} \\ =
 (b_2^{-1}-a_2)^{-1}(1-b^{-1}A_1(b_1)b^{-1}A_2(b_2))^{-1}b^{-1}
 (b_1^{-1}-a_1)^{-1}.  
\end{multline}
We want to use the Carl Neumann series in (\ref{inverti}) so we observe
that
\begin{eqnarray*}
b^{-1}A_i(b_i)& = & b^{-1}((b_i^{-1}-a_i)^{-1}-b) \\
              & = & (b_iE((1-a_ib_i)^{-1})))^{-1}b_i(1-a_ib_i)^{-1}-1 \\
              & = & \left(1+\sum_{n=1}^\infty
              E((a_ib_i)^n)\right)^{-1}\left(1 + \sum_{n=1}^\infty 
              (a_ib_i)^n\right) -1 \\
              & = & \left(1 + \sum_{k=1}^\infty \left(-\sum_{n=1}^\infty
              E((a_ib_i)^n)\right)^k\right)\left(1 + \sum_{n=1}^\infty 
              (a_ib_i)^n\right) -1 \\
              & = & \sum_{k=1}^\infty\left(-\sum_{n=1}^\infty
              E((a_ib_i)^n)\right)^k + \sum_{n=1}^\infty(a_ib_i)^n \\ &
              & +   \left(\sum_{k=1}^\infty\left(-\sum_{n=1}^\infty
              E((a_ib_i)^n)\right)^k\right)  \left(\sum_{n=1}^\infty(a_ib_i)^n\right).
\end{eqnarray*}
We infer that
\begin{equation}
\norm{b^{-1}A_i(b_i)} \leq \frac{\frac{2}{9}}{1-\tfrac{2}{9}}  +
\frac{\frac{2}{11}}{1-\frac{2}{11}} + 
\frac{\frac{2}{9}}{1-\tfrac{2}{9}} \frac{\frac{2}{11}}{1-\frac{2}{11}} <1.
\end{equation}
The Carl Neumann series now applies to (\ref{inverti}): 
\begin{multline} \label{Neumannseries}
 (b_3^{-1} -( a_1+a_2))^{-1} \\ =
 (b_2^{-1}-a_2)^{-1}\left( \sum_{n=0}^\infty
 \bigl(b^{-1}A_1(b_1)b^{-1}A_2(b_2)\bigr)^n \right) b^{-1}  
 (b_1^{-1}-a_1)^{-1}
\end{multline}
Since $A_1(b_1)$ and $A_2(b_2)$ are $\mathcal B$-free and centered 
(\ref{Neumannseries}) and (\ref{eq:1.4}) implies 
\begin{eqnarray*}
g_{a_1+a_2}(b_3) & = & E((b_3^{-1} -( a_1+a_2))^{-1}) \\ & = &
E((b_2^{-1}-a_2)^{-1}b^{-1}(b_1^{-1}-a_1)^{-1}) \\ 
 &=& E((b_2^{-1}-a_2)^{-1})b^{-1}E((b_1^{-1}-a_1)^{-1}) \\
& = & bb^{-1}b = b.
\end{eqnarray*} 

Thus by (\ref{b3def})
\begin{eqnarray*}
R_{a_1+a_2}(b) & = & (g_{a_1+a_2}^{<-1>}(b))^{-1} - b^{-1} \\
               & = & b_3^{-1}- b^{-1} \\
               & = & (b_1^{-1} - b^{-1}) + (b_2^{-1} - b^{-1}) \\
               & = & (g_{a_1}^{<-1>}(b))^{-1} -b^{-1}+
               (g_{a_2}^{<-1>}(b))^{-1} - b^{-1}\\     
               & = & R_{a_1}(b) + R_{a_2}(b).               
\end{eqnarray*}
\end{proof}

\section{Connection to Speicher's combinatorial approach}

Actually we do not have to restrict ourselves to invertible elements in
definition \ref{Rdef_inv}, because we will now show that the
$R$-transform has a $\mathcal B$-removable singularity in each
non-invertible element exactly as is the case for the scalar R-transform \cite[prop. 3.1]{Haa}.

For this we need the following description of the sum over all
irreducible non crossing partitions. We wil say that a partition $\pi\in \NC(r)$ is
irreducible if $1\sim_\pi r$ and we will denote concatenation of
non-crossing partitions by $\sqcup$.
\begin{lemma} \label{irrlemma}
let $\mathcal A$ be a unital Banach algebra, $1\in \mathcal
B\subset \mathcal A$ a unital Banach sub-algebra of $\mathcal A$, and 
$E:\mathcal A\to \mathcal B$ a conditional expectation. Then 
\begin{equation} \label{irreducible}
\sum_{j=1}^r (-1)^{j+1}\sum_{\substack{n_1+\cdots + n_j=r \\
n_1,\ldots,n_j \geq 1}} 
\sum_{\substack{\pi\in\NC(r)\\ \pi\leq 1_{n_1} \sqcup \cdots \sqcup
1_{n_j}}}  \kappa^\mathcal B_\pi((ba)^{\otimes r}) =
\sum_{\substack{\pi\in \NC(r) \\ 1\sim_\pi r}} \kappa^\mathcal
B_\pi((ba)^{\otimes r}),
\end{equation}
for all
$r\in \mathbb N$ , $a\in \mathcal A$ and $b\in \mathcal B$.
\end{lemma}

\begin{proof}
The proof is on induction in $r$. For $r=1$ (\ref{irreducible}) is
obvious since both sides of the equation reduces to $\kappa_1^\mathcal
B(ba)$. 

Now assume that (\ref{irreducible}) is true for all indices strictly
less that $r$. Then the left-hand side of (\ref{irreducible}) is
\begin{equation}\label{irr1}
E((ba)^r)+\sum_{j=2}^r (-1)^{j+1}\sum_{n_1=1}^{r-j+1}\sum_{\substack{n_2+\cdots
+ n_j=r-n_1 \\ 
n_2,\ldots,n_j \geq 1}} 
\sum_{\substack{\pi\in\NC(r)\\ \pi\leq 1_{n_1} \sqcup \cdots \sqcup
1_{n_j}}}  \kappa^\mathcal B_\pi((ba)^{\otimes r}).
\end{equation} 
Since ``$\sum_{j=2}^r \sum_{n_1=1}^{r-j+1}=
\sum_{n_1=1}^{r-1}\sum_{j=2}^{r-n_1+1}$'' we can rewrite (\ref{irr1}) as
\begin{multline*} 
E((ba)^r) + \left(\sum_{n_1=1}^{r-1} \sum_{\pi_1\in
\NC(n_1)}\kappa_{\pi_1}^\mathcal B((ba)^{\otimes n_1}) \sum_{j=2}^{r-n_1+1}
(-1)^{j+1}\right.
\\
\left.\sum_{\substack{n_2+\cdots 
+ n_j=r-n_1 \\ 
n_2,\ldots,n_j \geq 1}} 
 \sum_{\substack{\pi\in\NC(r-n_1)\\ \pi\leq 1_{n_2} \sqcup \cdots \sqcup
1_{n_j}}}  \kappa^\mathcal B_\pi((ba)^{\otimes (r-n_1)})\right).
\end{multline*}
Defining index $i:=j-1$ and $m_l:= n_{l+1}$ for $l\in\{1,\ldots,j-1\}$ we have
\begin{multline} \label{irr3}
E((ba)^r) + \left(\sum_{n_1=1}^{r-1} \sum_{\pi_1\in
\NC(n_1)}\kappa_{\pi_1}^\mathcal B((ba)^{\otimes n_1}) \right.
\\
\left.\left(- \sum_{i=1}^{r-n_1}
(-1)^{i+1}
\sum_{\substack{m_1+\cdots 
+ m_i=r-n_1 \\ 
m_1,\ldots,m_i \geq 1}} 
 \sum_{\substack{\pi\in\NC(r-n_1)\\ \pi\leq 1_{m_1} \sqcup \cdots \sqcup
1_{m_i}}}  \kappa^\mathcal B_\pi((ba)^{\otimes (r-n_1)})\right)\right).
\end{multline}
Using the induction hypotheses on the inner parenthesis of (\ref{irr3})
we have
\begin{multline} \label{irr4}
E((ba)^r) \\ - \left(\sum_{n_1=1}^{r-1} \sum_{\pi_1\in
\NC(n_1)}\kappa_{\pi_1}^\mathcal B((ba)^{\otimes n_1})
\left( \sum_{\substack{\pi\in \NC(n_1+1,\ldots,r) \\ (n_1+1)\sim_\pi
r}}  \kappa_\pi^{\mathcal B}((ba)^{\otimes (r-n_1)})\right)\right).
\end{multline}
But the (outer) parenthesis in (\ref{irr4}) is just the sum over all
non-irreducible 
partitions in $\NC(r)$, so (\ref{irr4}) reduces to
\begin{equation*}
\sum_{\substack{\pi\in \NC(r) \\ 1\sim_\pi r}} \kappa^\mathcal
B_\pi((ba)^{\otimes r}).
\end{equation*}
\end{proof}

We are now ready to remove the $\mathcal B$-singularities of the
$R$-transform. 

Let $g_a(b)\in \mathcal B(0,\tfrac{1}{11 \norm{a}})_{\inv}$ for some
$b\in \mathcal B(0,\tfrac{2}{11 \norm{a}})_{\inv}$. Then
\begin{eqnarray}
R_a(g_a(b)) & = & b^{-1}-(g_a(b))^{-1} \nonumber \\ 
            & = & b^{-1} - b^{-1}(E(1-ba)^{-1})^{-1} \nonumber \\
            & = & b^{-1} - b^{-1}\left(1 - \left(-\sum_{n=1}^\infty
            E((ba)^n)\right)\right)^{-1}  \nonumber \\
            & = &
            b^{-1}\sum_{j=1}^\infty(-1)^{j+1}\left(\sum_{n=1}^\infty  
            E((ba)^n)\right)^j \nonumber \\
            & = & \left(\sum_{n=0}^\infty E(a(ba)^n)\right) \left(1+
            \sum_{n=1}^\infty E((ba)^n)\right)^{-1} \nonumber \\ \label{Rdef}
            & = & E(a(1-ba)^{-1})(E(1-ba)^{-1})^{-1}.
\end{eqnarray}
We can thus adopt (\ref{Rdef}) as the definition of the $R$-transform
even for non-invertible $g_a(b)\in \mathcal
B(0,\tfrac{1}{11\norm{a}})$. Doing this we recover Speicher's definition
of the amalgamated $R$-transform from \cite[Th. 4.1.12]{AMS} by use of
lemma \ref{irrlemma}.
We have
\begin{eqnarray}
R_a(g_a(b)) & = &  \sum_{n=0}^\infty E(a(ba)^n)
            \left(1+\sum_{n=1}^\infty E((ba)^n)\right)^{-1} \nonumber \\ 
            & = &  \sum_{n=0}^\infty E(a(ba)^n) \left(\sum_{j=1}^\infty
            (-1)^{j+1} \left(\sum_{n=1}^\infty
            E((ba)^n)\right)^{j-1}\right)  \nonumber \\
            & = & b^{-1}\sum_{r=1}^\infty \left( \sum_{j=1}^r (-1)^{j+1}
            \sum_{\substack{n_1+ \cdots + n_j =r \\ n_1,\ldots,n_j\geq
            1}}E((ba)^{n_1})\cdots E((ba)^{n_j})\right) \label{R1}
\end{eqnarray} 
Actually the last line does not make sense for non-invertible elements,
but the singularity is obviously removable because $E$ has the $\mathcal B$-bi-module-property, and the following
computations becomes notationally simpler, when we write (\ref{R1}) in this way. We now use lemma \ref{irrlemma} in (\ref{R1})
\begin{eqnarray}
R_a(g_a(b)) & = &  b^{-1}\sum_{r=1}^\infty \left( \sum_{j=1}^r (-1)^{j+1}
            \sum_{\substack{n_1+ \cdots + n_j =r \\ n_1,\ldots,n_j\geq
            1}} \sum_{\substack{\pi\in \NC(r) \\ \pi\leq
            1_{n_1}\sqcup\cdots \sqcup 1_{n_j}}} \kappa_\pi^\mathcal
            B((ba)^{\otimes r}) \right) \nonumber \\
            & = &  b^{-1}\sum_{r=1}^\infty \left( \sum_{\substack{\pi\in
            \NC(r)\\ 1\sim_\pi r}} \kappa_\pi^\mathcal B((ba)^{\otimes
            r})\right)   \label{R2}          
\end{eqnarray} 
Rearranging the sum after the number of elements in the block that contains $1$ and $r$ (which are always in the same block), we actually just sum over all possible outer partitions. 
\begin{multline} \label{R2m}
  R_a(g_a(b))   =\\   \sum_{r=2}^\infty \sum_{j=2}^r \sum_{\substack{i_2+i_3+\cdots i_j +j =r \\ i_2, i_3, \ldots, i_j \geq 0}} \kappa_\pi^\mathcal B\left(a\otimes E((ba)^{i_2})ba \otimes E((ba)^{i_3})ba \otimes \cdots E((ba)^{i_j})ba\right) \\
= \sum_{k=2}^\infty\sum_{i_2, i_3, \ldots, i_k=0}^\infty \kappa_\pi^\mathcal B\left(a\otimes E((ba)^{i_2})ba \otimes E((ba)^{i_3})ba \otimes 
\cdots E((ba)^{i_k})ba\right).
\end{multline}
Rearranging terms in (\ref{R2m}) we get
\begin{eqnarray}
R_a(g_a(b)) & = & \sum_{r=1}^\infty \kappa_r^\mathcal B\left(a \otimes
\left(\left(\sum_{n=0}^\infty E((ba)^n)\right)ba\right)^{\otimes (r-1)}
\right) \nonumber \\
 & = & \sum_{r=1}^\infty \kappa_r^\mathcal B(a \otimes g_a(b)a\otimes \cdots
 \otimes g_a(b)a). \label{R3}
\end{eqnarray}
This is exactly Speicher's way of defining the amalgamated
$R$-transform, and additivity of the $R$-transform on $\mathcal B$-free
variables 
thus follows from \cite[th. 4.1.7]{AMS}. The new thing is that we have
produced a concrete neighborhood where the sum in (\ref{R3}) makes
sense, that is, if $g_a(b)\in \mathcal 
B(0,\tfrac{1}{11\norm{a}})$ for some $b\in \mathcal
B(0,\tfrac{2}{11\norm{a}})$ then (\ref{R3}) is convergent.

\section{Amalgamated S-transform in Banach-algebras}

Again we let $\mathcal A$ be a unital Banach algebra, 
$1\in \mathcal B\subset \mathcal A$ be a unital Banach sub-algebra of
$\mathcal A$ and $E:\mathcal A\to \mathcal B$ a conditional expectation. The results on the amalgamated $S$-transform are obtained via a similar approach as in the last section.

Let $a\in \mathcal A$ be a fixed
element. Define $\Psi_a :\mathcal B(0,\tfrac{1}{\norm{a}}) \to
\mathcal B$ by
\begin{gather}
 \Psi_a(b) = \sum_{n=1}^\infty E((ba)^n) = E((1-ba)^{-1})-1,\label{2.2}
\end{gather}
for $b\in \mathcal B(0,\tfrac{1}{\norm{a}})$.

We proceed as in the previous section to show that  $\Psi_a$
is injective in a neighborhood $\mathcal
B(0,\tfrac{\alpha}{\norm{a}})$ and maps $\mathcal
B(0,\tfrac{\alpha}{\norm{a}})$ onto a neighborhood of $0$ which
contains a neighborhood $\mathcal
B(0,\tfrac{\beta}{\norm{a}})$, where $\alpha$ and $\beta$ are constants
to be determined.

\begin{lemma}
Let $a\in\mathcal A$ such that $E(a)\in \mathcal B_{\inv}$ and let
$\Psi_a:\mathcal 
B(0,\tfrac{1}{\norm{a}})\to \mathcal B$ be the 
function defined by (\ref{2.2}). Then
$\Psi_a$ is 1-1 on $\mathcal 
B(0,\tfrac{1}{4\norm{a}^2\norm{E(a)^{-1}}})$. 
\end{lemma}
 
\begin{proof} 
Define $\Gamma_a :\mathcal B(0,\tfrac{1}{\norm{a}})\to \mathcal B$ by 
$\Gamma_a : b \mapsto \Psi_a(b) E(a)^{-1}$. 
Let $b_1,b_2\in \mathcal B(0,\tfrac{\alpha}{\norm{a}^2\norm{E(a)^{-1}}})$ where
$0<\alpha<1$ is to be determined. Observe that since $1\leq
\norm{a}\norm{E(a)^{-1}}$ we have $ \mathcal
B(0,\tfrac{\alpha}{\norm{a}^2\norm{E(a)^{-1}}}) \subseteq \mathcal
B(0,\tfrac{\alpha}{\norm{a}})$. Now
\begin{eqnarray*}
 \Gamma_a(b_1) -\Gamma_a(b_2) & = & b_1-b_2 + \bigl( E(b_1ab_1a)-E(b_2ab_2a) \\ 
   & & +
 E(b_1ab_1ab_1a)- E(b_2ab_2ab_2a) + \cdots)\bigr) E(a)^{-1} \\
 & = & b_1-b_2 + \bigl(E((b_1-b_2)ab_1a) + E(b_2a(b_1-b_2)a) \\
 &   & + E(b_1-b_2)ab_1ab_1a) + E(b_2a(b_1-b_2)ab_2a)  \\ &&+
   E(b_2ab_2a(b_1-b_2)a) + \cdots\bigr)E(a)^{-1}.
\end{eqnarray*}
Defining $c= \norm{a}\max\{\norm{b_1},\norm{b_2}\}$ we estimate
\begin{multline*}
\norm{\Psi_a(b_1)-\Psi_a(b_2) -(b_1-b_2)}  \\ \leq 
2\norm{a}\max\{\norm{b_1},\norm{b_2}\}\norm{b_1-b_2}\norm{a}\norm{E(a)^{-1}} \\
 + 3
 \norm{a}^2\max\{\norm{b_1},\norm{b_2}\}^2\norm{b_1-b_2}\norm{a}\norm{E(a)^{-1}} \\  
  + \cdots  \\
 \leq  \norm{b_1-b_2}
 \norm{a}\norm{E(a)^{-1}}\left(\sum_{n=0}^\infty(n+1)c^n -1\right) 
\\
 =
 \left(\frac{2-c}{(1-c)^2}\right)c\norm{a}\norm{E(a)^{-1}}\norm{b_1-b_2}
 \\
< \frac{\alpha(2-\alpha)}{(1-\alpha)^2}\norm{b_1-b_2}.
\end{multline*}
We infer that when $\alpha < \tfrac{1}{4}$ then
$\tfrac{\alpha(2-\alpha)}{(1-\alpha)^2}< 1$ and the lemma follows. 
\end{proof}

\begin{lemma} \label{differentiale2} Let $a\in \mathcal A$ be fixed.
  Let $\Psi_a:\mathcal B(0,\tfrac{1}{\norm{a}})\to \mathcal B$ be
  the function defined by (\ref{2.2}). Then $\Psi_a$ is
  Fr\'echet differentiable, and the differential  
of $\Psi_a$ is 
\begin{gather*}
  D\Psi_a(b):h\mapsto E((1-ba)^{-1}ha(1-ba)^{-1})
\end{gather*}
for $b\in \mathcal B(0,\tfrac{1}{\norm{a}})$ and $h\in \mathcal B$.
Furthermore the differential, $D\Psi_a:\mathcal
B(0,\tfrac{1}{\norm{a}})\to  
\mathcal B$, is continuous, so that $\Psi_a$ is
differentiable of class  
$\C^1$. If $b \in \mathcal
B(0,\tfrac{1}{\norm{a}})$ then
\begin{gather}
  \label{DPsiavurdering}
  \norm{D\Psi_a(b)-D\Psi_a(0)} < \frac{\norm{b}\norm{a}^2(2-\norm{b}\norm{a})}{(1-\norm{b}\norm{a})^2}. 
\end{gather} 
\end{lemma}

\begin{proof}
Recall that for $b_1,b_2\in
\mathcal B(0,\tfrac{1}{\norm{a}})$ we have
\begin{equation*}
(1-b_1a)^{-1}-(1-b_2a)^{-1} = (1-b_1a)^{-1}(b_1-b_2)a(1-b_2a)^{-1},
\end{equation*}
so 
\begin{multline} \label{Psialemma1}
\norm{(1-b_1a)^{-1}-(1-b_2a)^{-1}}  \\ \leq
\frac{1}{1-\norm{b_1}\norm{a}}\frac{1}{1-\norm{b_2}\norm{a}}
\norm{a}\norm{b_1-b_2} \to 0
\end{multline}
for $b_1\to b_2$ in norm. Also
\begin{eqnarray}
\Psi_a(b_1)-\Psi_a(b_2) & = & \sum_{n=1}^\infty E((b_1a)^n) -
\sum_{n=1}^\infty E((b_2a)^n) \nonumber \\
                        & = & E\left((1-b_1a)^{-1} - (1-b_2a)^{-1}\right) \nonumber
                        \\
                        & = &
                        E\left((1-b_1a)^{-1}(b_1-b_2)a(1-b_2a)^{-1}\right). 
                        \label{Psialemma2}  
\end{eqnarray} 
Combining (\ref{Psialemma1}) and (\ref{Psialemma2}) we get
\begin{multline*}
\Psi_a(b+h)-\Psi_a(b)  =   E\left((1-(b+h)a)^{-1}ha(1-ba)^{-1}\right) \\
             =   E\left((1-ba)^{-1}ha(1-ba)^{-1}\right) \\
                +  E\left(\left((1-(b+h)a)^{-1}- (1-ba)^{-1}\right)ha(1-ba)^{-1}\right) \\
             =    E\left((1-ba)^{-1}ha(1-ba)^{-1}\right) + O(\norm{h}^2)
\end{multline*}
for all $b \in \mathcal B(0,\tfrac{1}{\norm{a}})$ and $h\in \mathcal B$
of small norm. This shows the first part of the lemma.

Let $b_1,b_2\in \mathcal B(0,\tfrac{1}{\norm{a}})$. Continuity of
$D\Psi_a:b\mapsto D\Psi_a(b)$ follows from
\begin{multline} \label{Dpsia}
\norm{D\Psi_a(b_1)-D\Psi_a(b_2)} \\ =  \sup_{\norm{h}\leq
1}\norm{D\Psi_a(b_1)(h)-D\Psi_a(b_2)(h)} \\
      = \sup_{\norm{h}\leq 1} \norm{E\left((1-b_1a)^{-1}ha(1-b_1a)^{-1}\right)-
      E\left((1-b_2a)^{-1}ha(1-b_2a)^{-1}\right) } \\
      \leq \sup_{\norm{h}\leq 1}
      \norm{(1-b_1a)^{-1}ha(1-b_1a)^{-1}- (1-b_1a)^{-1}ha(1-b_2a)^{-1}}
      \\
      +  \leq \sup_{\norm{h}\leq 1}
      \norm{(1-b_1a)^{-1}ha(1-b_2a)^{-1}- (1-b_2a)^{-1}ha(1-b_2a)^{-1}}
      \\
    \leq
    \norm{(1-b_1a)^{-1}-(1-b_2a)^{-1}}\norm{a}\left(\norm{(1-b_1a)^{-1}}+
    \norm{(1-b_2a)^{-1}}\right) \\ \to 0
\end{multline}
for $b_1 \to b_2$ in norm. Specifically letting $b_2=0$ and $b_1=b$ in
(\ref{Dpsia}) we have
\begin{multline*}
\norm{D\Psi_a(b)-D\Psi_a(0)}  \leq \norm{(1-ba)^{-1}-1}\norm{a}
\left(\norm{(1-ba)^{-1}}+1\right) \\ 
     \leq  \norm{a} \frac{\norm{b}\norm{a}}{1-\norm{b}\norm{a}}\left(
      \frac{1}{1-\norm{b}\norm{a}} +1 \right) \\
     = 
      \frac{\norm{b}\norm{a}^2(2-\norm{b}\norm{a})}{
      (1-\norm{b}\norm{a})^2}.  
\end{multline*}
\end{proof}

\begin{proposition} \label{domain2}
Let $\mathcal A$ be a unital Banach algebra, $1\in \mathcal
B\subset \mathcal A$ a unital Banach sub-algebra and $E:\mathcal A\to
\mathcal B$ a
conditional expectation. Let $a\in \mathcal A$ be a fixed element and
assume that $E(a)\in \mathcal B_{\inv}$. The
function $\Psi_a:\mathcal 
B(0,\tfrac{1}{\norm{a}})\to \mathcal B$ defined by
$\Psi_a(b)=E((1-ba)^{-1})-1$ for $b\in
\mathcal B(0, 
\tfrac{1}{\norm{a}})$ is a
bijection of the neighborhood $\mathcal
B(0,\tfrac{1}{4\norm{a}^2\norm{E(a)^{-1}}})$ 
onto a neighborhood of $0$ which contains $\mathcal
B(0,\tfrac{1}{11\norm{a}^2\norm{E(a)^{-1}}^2})$. Furthermore   
\begin{gather}
\Psi_a^{<-1>}\left(\mathcal
   B(0,\tfrac{1}{11\norm{a}^2\norm{E(a)^{-1}}^2}) \right) 
   \subseteq \mathcal B(0,\tfrac{2}{11\norm{a}^2\norm{E(a)^{-1}}}) \\
\Psi_a^{<-1>}\left(\mathcal
   B(0,\tfrac{1}{11\norm{a}^2\norm{E(a)^{-1}}^2})_{\inv} \right) 
   \subseteq \mathcal B(0,\tfrac{2}{11\norm{a}^2\norm{E(a)^{-1}}})_{\inv} 
   \label{neighboor2} 
\end{gather}
\end{proposition}

\begin{proof}
 The proof is very similar to the proof of proposition \ref{domain} and is
 again an inspection of some of the proof of the Inverse Function
 theorem. We only give the changes.

Define $\Gamma_a:\mathcal B(0,\tfrac{1}{\norm{a}})\to \mathcal B$ by
$\Gamma_a:b\mapsto \Psi_a(b)E(a)^{-1}$. Then we have $D\Gamma(0)(h) =
h$. We now define $T:\mathcal B(0,\tfrac{1}{\norm{a}}) \to \mathcal B$ by
$T(b) = b - 
\Gamma_a(b)$ for $b\in \mathcal B(0,\tfrac{1}{\norm{a}})$ and observe
that $T(0)=0$ and $DT(0)=0$. Thus by use of (\ref{DPsiavurdering}) we have
\begin{eqnarray*}
 \norm{DT(0)} & = & \norm{D\Gamma_a(b)-D\Gamma_a(0)} \\
              & \leq & \norm{E(a)^{-1}}\norm{D\Psi_a(b)-D\Psi_a(0)} \\
              & \leq &
              \norm{E(a)^{-1}}\norm{a}^2\norm{b}\frac{2-\norm{b}\norm{a}}
              {(1-\norm{b}\norm{a})^2} < \frac{40}{81} <\frac 1 2 
\end{eqnarray*} 
for $b\in \mathcal B(0,\tfrac{2}{11\norm{a}^2\norm{E(a)^{-1}}})$.

We can now proceed exactly as in the proof of proposition \ref{domain}
to show that $\Gamma_a$ maps $\mathcal
B(0,\tfrac{2}{11\norm{a}^2\norm{E(a)^{-1}}})$ injectively onto a neighborhood of $0$
containing $\mathcal B(0,\tfrac{1}{11\norm{a}^2\norm{E(a)^{-1}}})$. Thus
if $b_0\in \mathcal B(0,\tfrac{1}{11\norm{a}^2\norm{E(a)^{-1}}^2})$ it is
obvious that $b_0E(a)^{-1} \in \mathcal
B(0,\tfrac{1}{11\norm{a}^2\norm{E(a)^{-1}}})$ so there exists $b\in
\mathcal B(0,\tfrac{2}{11\norm{a}^2\norm{E(a)^{-1}}})$ such that
$\Gamma_a(b)=b_0E(a)^{-1}$. But then $\Psi_a(b) = \Gamma_a(b)E(a)= b_0$, so $\Psi_a$ maps
$\mathcal B(0,\tfrac{2}{11\norm{a}^2\norm{E(a)^{-1}}})$ injectively onto a
neighborhood of $0$ containing $\mathcal
B(0,\tfrac{1}{11\norm{a}^2\norm{E(a)^{-1}}^2})$.

To see (\ref{neighboor2}) observe that for $b\in \mathcal
B(0,\tfrac{2}{11\norm{a}^2\norm{E(a)^{-1}}})$ we have  
\begin{eqnarray*}
\norm{E(a)^{-1}\sum_{n=1}^\infty E(a(ba)^n)} & \leq &  \norm{E(a)^{-1}}
\norm{a}^2 \norm{b} \sum_{n=0}^\infty (\norm{b}\norm{a})^n \\
& < &  \frac{\frac{2}{11}}{1-\frac{2}{11}} <1, 
\end{eqnarray*}
so (\ref{neighboor2}) now follows from 
\begin{equation} \label{Sinvvurdering} 
\Psi_a(b) = \sum_{n=1}^\infty E((ba)^n) = bE(a)\left(1 +
E(a)^{-1}\sum_{n=1}^\infty E(a(ba)^n)\right).
\end{equation}
\end{proof}

We can now define the amalgamated $S$-transform 
\begin{definition} \label{Sdef_inv}
  Let $\mathcal A$ be a unital Banach algebra, $1\in \mathcal
  B\subset \mathcal A$ a unital commutative Banach sub-algebra and 
  $E:\mathcal A\to \mathcal B$ a conditional expectation. Let $a\in \mathcal 
  A$ be a fixed non-zero element such that $E(a)\in \mathcal
  B_{\inv}$. Define the amalgamated $S$-transform of $a$ by 
  \begin{equation}
    \label{S-transform}
    S_a(b) := b^{-1}(1+b) \Psi_a^{<-1>}(b)
  \end{equation}
for $b\in \mathcal B(0, \tfrac{1}{11\norm{a}^2\norm{E(a)^{-1}}^2})_{\inv}$.
\end{definition}

We have the following amalgamated version of (\ref{Stransformen}).

\begin{theorem}
 Let $\mathcal A$ be a unital Banach algebra, $1\in \mathcal
  B\subset \mathcal A$ a unital commutative Banach sub-algebra and 
  $E:\mathcal A\to \mathcal B$ a conditional expectation. Let
  $a_1,a_2\in \mathcal  A$ be $\mathcal B$-free fixed non-zero elements
  such that 
  $E(a_1),E(a_2) \in \mathcal B_{\inv}$.
Let \\
$b\in \mathcal B(0,
  \min(\tfrac{1}{11\norm{a_1}^2\norm{E(a_1)^{-1}}^2}, 
 \tfrac{1}{11\norm{a_2}^2\norm{E(a_2)^{-1}}^2}, 
  \tfrac{1}{11\norm{a_1+a_2}^2\norm{E(a_1+a_2)^{-1}}^2}))_{\inv}$. Then
  \begin{equation*}
  S_{a_1 a_2}(b) = S_{a_1}(b) S_{a_2}(b).
  \end{equation*}
\end{theorem}

\begin{proof}
Let \\ $b\in \mathcal B(0,
  \min(\tfrac{1}{11\norm{a_1}^2\norm{E(a_1)^{-1}}^2}, 
 \tfrac{1}{11\norm{a_2}^2\norm{E(a_2)^{-1}}^2}, 
  \tfrac{1}{11\norm{a_1+a_2}^2\norm{E(a_1+a_2)^{-1}}^2}))_{\inv}$, and let 
 $b_1\in\mathcal B(0,\tfrac{2}{11\norm{a_1}^2\norm{E(a_1)^{-1}}})_{\inv}$ and $b_2 \in\mathcal B(0,\tfrac{2}{11\norm{a_2}^2\norm{E(a_2)^{-1}}})_{\inv}$ be the uniquely determined elements such that
  \begin{equation*}
    b = \Psi_{a_1}(b_1) = \Psi_{a_2}(b_2).
  \end{equation*}
Note that 
\begin{equation} \label{bomskriv}
  b+1 = E((1-b_1a_1)^{-1}) = E((1-a_2b_2)^{-1}),
\end{equation}
where the last equality follows since $\mathcal B$ is commutative.
Define
\begin{eqnarray*}
  A_1(b_1)& = &  (1-b_1a_1)^{-1} - E((1-b_1a_1)^{-1}), \\
  A_2(b_2)& = &  (1-a_2b_2)^{-1} - E((1-a_2b_2)^{-1}).
\end{eqnarray*}
By (\ref{bomskriv}) we have
\begin{eqnarray*}
  (1-b_1a_1)A_1(b_1) & = & 1 - (1-b_1a_1)(1+b) \\
  A_2(b_2)(1-a_2b_2) & = & 1 - (1+b)(1-a_2b_2),
\end{eqnarray*}
and thus
\begin{multline*}
  (1-b_1a_1)b\left(1-b^{-1}A_1(b_1)(1+b)^{-1}A_2(b_2)\right)(1-a_2b_2) \\
 = (1-b_1a_1)b(1-a_2b_2) - \\
 \bigl(1-(1-b_1a_1)(1+b)\bigr)(1+b)^{-1}\bigl(1-(1+b)(1-a_2b_2)\bigr) \\  
 = -(1+b)^{-1} + (1-b_1a_1) + (1-a_2b_2)  -(1-b_1a_1)(1-a_2b_2) \\
 = 1- (1+b)^{-1} - b_1a_1a_2b_2=\frac{b}{1+b}\left(1-\frac{1+b}{b}b_1a_1a_2b_2\right),
\end{multline*}
where $\frac{b}{1-b}$ and $\frac{1+b}{b}$ denotes the elements $(1+b)^{-1} b$ and $b^{-1}(1+b)$ respectively.
We claim that $\norm{b^{-1}A_1(b_1)(1+b)^{-1}A_2(b_2)}<1$. Inverting we have
\begin{multline*}
  \left(1- \frac{1+b}{b}b_1a_1a_2b_2\right)^{-1}\frac{1+b}{b} \\
 =  (1-a_2b_2)^{-1}\left(1-b^{-1}A_1(b_1)(1+b)^{-1}A_2(b_2)\right)^{-1}b^{-1}(1-b_1a_1)^{-1} \\
(1-a_2 b_2)^{-1} \sum_{n=0}^\infty ( b^{-1} A_1(b_1) (1+b)^{-1} A_2(b_2))^n b^{-1} (1-b_1 a_1)^{-1}.
\end{multline*}
Since $A_1(b_1)$ and $A_2(b_2)$ are $\mathcal B$-free and $E(A_1(b_1))=
E(A_2(b_2))=0$, $\mathcal B$-freeness implies by (\ref{eq:1.4}) that
\begin{multline*}
  E((1-\frac{1+b}{b} b_1 a_1 a_2 b_2)^{-1})\frac{1+b}{b} \\ = E((1-a_2
  b_2)^{-1})b^{-1}  
E((1-b_1a_1)^{-1}) = \frac{(1+b)^2}{b},
\end{multline*}
so we conclude that
\begin{equation*}
  \Psi_{a_1a_2}(\frac{1+b}{b}b_1b_2) = b.
\end{equation*}
Hence
\begin{equation*}
  S_{a_1a_2}(b) = \frac{1+b}{b}\left(\frac{1+b}{b}b_1b_2\right) = 
\left(\frac{1+b}{b}b_1\right)\left(\frac{1+b}{b}b_2\right) = S_{a_1}(b)S_{a_2}(b).
\end{equation*}
To prove the claim assume for a moment that 
\begin{equation}\label{1mindreend2}
\norm{a_1}\norm{E(a_1)^{-1}} \leq \norm{a_2}\norm{E(a_2)^{-1}} 
\end{equation}
and observe that
\begin{multline*}
  b^{-1}A_1(b_1)   =  \Psi_{a_1}(b_1)^{-1}\left(\sum_{n=1}^\infty (b_1a_1)^n 
 - \Psi_{a_1}(b_1) \right) \\
=  \left(b_1E(a_1)\left(1+E(a_1)^{-1}\sum_{k=1}^\infty E(a_1(b_1a_1)^k)\right)\right)^{-1}\sum_{n=1}^\infty (b_1a_1)^n - 1 \\
 = \left(1+E(a_1)^{-1}\sum_{k=1}^\infty E(a_1(b_1a_1)^k)\right)^{-1}E(a_1)^{-1}a_1 \sum_{n=0}^\infty (b_1a_1)^n -1,
\end{multline*}
so
\begin{multline}
  \label{vurd1}
  \norm{b^{-1}A_1(b_1)}< \left(\frac{1}{1-\frac{2}{9}}\right) \norm{a_1}\norm{E(a_1)^{-1}}\left(\frac{1}{1-\frac{2}{11}}\right) +1 \\
= \frac{11}{7} \norm{a_1}\norm{E(a_1)^{-1}}+1 \leq \frac{11}{7} \norm{a_2}\norm{E(a_2)^{-1}}+1.
\end{multline}
Also
\begin{eqnarray*}
  (1+b)^{-1}A_2(b_2) & = & (1+b)^{-1}(1-a_2b_2)^{-1} -1 \\
   & = & \sum_{n=0}^\infty (-b)^n\sum_{k=0}^\infty (a_2b_2)^k -1 \\
   & = & \sum_{n=1}^\infty(-b)^n + \sum_ {k=1}^\infty (a_2b_2)^k + \sum_{n=1}^\infty(-b)^n\sum_ {k=1}^\infty (a_2b_2)^k,
\end{eqnarray*}
so
\begin{multline} \label{vurd2}
  \norm{ (1+b)^{-1}A_2(b_2)} \\ < \left(\frac{1}{10} + \frac{2}{9} + 
\frac{1}{10}\frac{2}{9}\right)\frac{1}{\norm{a_2}\norm{E(a_2)^{-1}}} = \frac{32}{90} \frac{1}{\norm{a_2}\norm{E(a_2)^{-1}}}.
\end{multline}
The claim now follows by combining (\ref{vurd1}) and (\ref{vurd2}) 
\begin{multline*}
 \norm{b^{-1}A_1(b_1)(1+b)^{-1}A_2(b_2)} \\ \leq  \norm{b^{-1}A_1(b_1)}  \norm{ (1+b)^{-1}A_2(b_2)} \\
< \left(\frac{11}{7}\norm{a_2}\norm{E(a_2)^{-1}}+1\right) \frac{32}{90}\frac{1}{\norm{a_2}\norm{E(a_2)^{-1}}} \\ <\frac{11}{7}\frac{32}{90}+\frac{32}{90}<1.
\end{multline*}
Finally if we have $\norm{a_1}\norm{E(a_1)^{-1}} \leq \norm{a_2}\norm{E(a_2)^{-1}}$ we just do similar calculations on $A_1(b_1)(1+b)^{-1}A_2(b_2)b^{-1}$ instead. 
\end{proof}

\section{Examples}

%\begin{example} \label{counterSmultiplicative}
%There exists counterexample to multiplicativity of $S$-transform in the
%case of a non-commutiatve amalgamated algebra. One can be found  in
%$3\times 3$-matrices with conditional expectation given by
%\begin{equation*}
%  E:\begin{bmatrix}
%  a_{11} & a_{12} & a_{13} \\
%  a_{21} & a_{22} & a_{23} \\
%  a_{31} & a_{32} & a_{33} 
%  \end{bmatrix} \mapsto 
%  \begin{bmatrix}
%   a_{11} & 0 & a_{13} \\
%  0 & a_{22} & 0 \\
%  a_{31} & 0 & a_{33}.
%  \end{bmatrix}
%\end{equation*}
%Details to be filled in.........
%\end{example}

The amalgamated $R$- and $S$-transform is related as follows. 
\begin{example}
 Let $\mathcal A$ be a unital Banach algebra, and let $1\in \mathcal
  B\subset \mathcal A$ be a unital commutative Banach sub-algebra, and 
  $E:\mathcal A\to \mathcal B$ a conditional expectation. Let $a\in \mathcal 
  A$ be a fixed non-zero element such that $E(a)\in \mathcal
  B_{\inv}$. As in the scalar case \cite{HL,NS2} we have the following relation beween
  the amalgamated $R$- and $S$-transform 
  \begin{equation}
  \label{RSrelation}
  bS_a(b) = [bR_a(b)]^{<-1>}.
  \end{equation}
To see (\ref{RSrelation}) note that for $b\in \mathcal B_{\inv}$ of
small norm we have
\begin{eqnarray*}
   g_a(b)R_a(g_a(b))& = & g_a(b)(b^{-1}-g_a(b)^{-1}) \\
                    & = & E((1-ba)^{-1})bb^{-1}-1 \\
                    & = & \sum_{n=1}^\infty E((ba)^n) = \Psi_a(b).
\end{eqnarray*} 
and 
\begin{eqnarray*}
\Psi_a(b)S_a(\Psi_a(b)) & = & \Psi_a(b)\frac{1+\Psi_a(b)}{\Psi_a(b)}b \\
                        & = & \left(1+\sum_{n=1}^\infty
                        E((ba)^n)\right)b \\
                        & = & g_a(b).
\end{eqnarray*}
So (\ref{RSrelation}) now follows easily because
\begin{equation*}
g_a(b)R_a(g_a(b))S_a\bigl(g_a(b)R_a(g_a(b))\bigr ) =
\Psi_a(b)S_a(\Psi(b)) = g_a(b)  
\end{equation*}
and
\begin{equation*}
\Psi_a(b)S_a(\Psi_a(b))R_a\bigl(\Psi_a(b)S_a(\Psi_a(b))\bigr)=
g_a(b)R_a(g_a(b))= 
\Psi_a(b). 
\end{equation*}
\end{example}
\begin{flushright}
  \qedsymbol
\end{flushright}

We have the following dilation formula.

\begin{example}[Dilations]
 Let $\mathcal A$ be a unital Banach algebra, and let $1\in \mathcal
 B\subset \mathcal A$ be a unital commutative Banach sub-algebra, and 
  $E:\mathcal A\to \mathcal B$ a conditional expectation. Let $a\in \mathcal 
  A$ be a fixed non-zero element such that $E(a)\in \mathcal
  B_{\inv}$. Assume that $z\in \mathcal B_{\text{inv}}$.

 Then for $b\in \mathcal B(0,\tfrac{1}{\norm{za}})$ and $bz\in \mathcal B(0,\tfrac{1}{\norm{a}})$  we have
  \begin{equation*}
  \Psi_{za}(b) = \sum_{n=1}^\infty E((b(za))^n) = \sum_{n-1}^\infty E(((bz)a)^n) = \Psi_a(bz),
  \end{equation*} 
and
\begin{equation*}
  S_z(\Psi_z(b)) = \left(\sum_{n=1}^\infty (bz)^n \right)^{-1}\left(\sum_{n=0}^\infty (bz)^n\right)bzz^{-1}=z^{-1},
\end{equation*}
for $\norm{b}$ sufficiently small. Thus
\begin{eqnarray*}
S_{za}(\Psi_{za}(b)) &=& \Psi_{za}(b)^{-1}(1+\Psi_{za}(b))  b 
\\ &=& \Psi_{za}(b)^{-1}(1+\Psi_{za}(b))(bz)z^{-1} \\ 
&=& S_a(\Psi_{a}(bz))S_z(\Psi_a(bz)) \\ &=& S_a(\Psi_{za}(b))S_z(\Psi_{za}(b))
\end{eqnarray*}
so 
\begin{equation*}
  S_{za}(b)= S_a(b)S_z(b)
\end{equation*}
for $b$ invertible and $\norm{b}$ sufficiently small. 
\end{example}
\begin{flushright}
  \qedsymbol
\end{flushright}
Note that in the above example we did \emph{not} use the commutativeness of $\mathcal B$, so the example above shows that if we define the amalgamated $S$-transform by (\ref{S-transform}) for non-commutative $\mathcal B$ we would actually have to look for a product formula of the form:
\begin{equation} \label{non-com}
  S_{a_1a_2}(b) = S_{a_2}(b) S_{a_1}(b)
\end{equation}
for $a_1$ $*$-free from $a_2$ and $\norm{b}$ sufficiently small. This suggests that it is actually more natural to consider the inverse of the $S$-transform than the $S$-transform.

Unfortunately we do not know whether (\ref{non-com}) is true or not when $\mathcal B$ is non-commutative, but our guess is that it is not true in general.

\providecommand{\bysame}{\leavevmode\hbox to3em{\hrulefill}\thinspace}


\begin{thebibliography}{25}
 

%\bibitem[DH1]{DT} K. Dykema, U. Haagerup, \emph{DT-operators and
%decomposability of Voiculescu's circular operator}, preprint (2002).

%\bibitem[DH2]{invsub} \bysame, \emph{Invariant subspaces of the
%quasinilpotent DT-operator}, preprint (200?).


%\bibitem[HP]{HiaiPetz} F. Hiai, D. Petz, \emph{The Semicircle Law, Free
%random Variables and Entropy}, Mathematical Surveys and Monographs,
%Vol. 77, Am Math. Soc., (2000).

\bibitem[Haa]{Haa} U. Haagerup, \emph{On Voiculescu's $R$- and
$S$-transform for free non-commuting random variables}, Fields Institute
Communications, Vol. \textbf{12} , AMS, p. 127-148 (1997).

\bibitem[HL]{HL} U. Haagerup, F. Larsen, \emph{Brown's Spectral Distribution Measure for R-diagonal Elements in Finite von Neumann Algebras}, J. Funct. Anal. \textbf{176}, 
331-367 (2000).

\bibitem[VLH]{VLH} V. L. Hansen, \emph{Fundamental Concepts in Modern Analysis}, World Scientific, 1999.

\bibitem[NS1]{NS1} A. Nica, R. Speicher, \emph{A ``Fourier transform'' for multiplicative functions on non-crossing partitions}, J. of Algebraic Combinatorics \textbf{6} (1997), 141-160.

\bibitem[NS2]{NS2} \bysame, \emph{R-diagonal elements -a common approach to Haar 
unitaries and circular elements}, Fields Institute Communications, Vol. \textbf{12},
 AMS, p. 149-188 (1997).

\bibitem[Sp1]{Sp1} R. Speicher, \emph{Multiplicative functions on the lattice of non-crossing partitions and free convolution}, Math. Annalen., \textbf{298}, 611-628, 1994.
\bibitem[Sp2]{AMS} \bysame, \emph{Combinatorial theory of the free product with amalgamation and operator-valued free probability theory}, Mem. Amer. Math. Soc. \textbf{132} (1998), no 627, x+88 pp. 

\bibitem[Sp3]{Combinatorics} \bysame, \emph{Combinatorics of free probability
theory}, ``Free probability and operator spaces'', IHP, Paris, 1999.


\bibitem[Voi1]{Voi1} D. Voiculescu, \emph{Addition of certain non-commuting random variables}, Journal of Functional Analysis, \textbf{66}, 323-346, 1986.

\bibitem[Voi2]{Voi2} \bysame, \emph{Multiplication of certain non-commuting random variables}, Journal of Operator Theory, \textbf{18}, 223-346, 1987.

\bibitem[Voi3]{Voi3} \bysame, \emph{Operations on certain non-commuting operator-valued random variables}, Ast\'erisque, \textbf{232}, 243-275, 1995.

%\bibitem[Sn1]{SniadyDT} P. \'Sniady, \emph{Inequality for Voiculescu's
%entropy in terms of Brown measure},
%Internat. Math. Res. Notices \textbf{2003}, 51-64. 

%\bibitem[Sn2]{Sniady} \bysame, \emph{Multinomal identities arising from
%the free probability}, J. Comb. Th. A, no. 1, \textbf{101} (2003), 1-19.

%\bibitem[SnSp]{SnSp} P. \'Sniady, R. Speicher, \emph{Continuous family
%of invarinat subspaces for R-diagonal operators},
%Invent. Math. \textbf{146} (2001), no. 2, 329-363.

%\bibitem[NSS1]{NSS1} A. Nica, D. Shlyakthenko, R. Speicher, \emph{Some
%minimization problems for the free analogue of the free Fischer
%information}, Ad.Math. \textbf{141} (1999), no. 2, 282-321.

%\bibitem[NSS2]{NSS2} \bysame, \emph{Operator-valued
%distributions. I. Characterizations of freeness},
%Int. Math. Res. Not. \textbf{2002}, no. 29, 1509-1538.

%\bibitem[Voi1]{Voi1} D. Voiculescu, \emph{The analouges of entropy and
%of Fischer's information measure in free probability theory, I}
%Communications Math. Physics \textbf{155} (1993), 71-92.

%\bibitem[Voi2]{Voi2} \bysame,  \emph{The analouges of entropy and
%of Fischer's information measure in free probability theory, II},
%Invent. Math. \textbf{118} (1994) 411-440.

%\bibitem[Voi5]{Voi5} \bysame, \emph{The analouges of entropy and
%of Fischer's information measure in free probability theory, V:
%Non-commutative Hilbert Transforms}, Inventiones Math. \textbf{132}
%(1998), 189-227.

%\bibitem[Voi6]{Voi6} \bysame, \emph{The analouges of entropy and
%of Fischer's information measure in free probability theory, VI:
%Liberation and Mutual Free Information}, Advances in Mathematics
%\textbf{146} (1999), 101-166.

\end{thebibliography}
\end{document}